\documentclass[12pt]{article}
\usepackage{amssymb}
\usepackage{amsmath,amsthm}
\usepackage[latin1]{inputenc}
\usepackage{citesort}
\usepackage{graphicx}
\usepackage{epsfig}

 \newtheorem{remark}{Remark}

 \newtheorem{theorem}[remark]{Theorem}

  \newtheorem{claim}[remark]{Claim}

\title{Spectral study of alliances in graphs}

\author{J. A. Rodr\'{\i}guez\footnote{e-mail:\mbox{\tt
    juanalberto.rodriguez\@@urv.net}} \\
{\em Department of Computer Engineering and Mathematics}\\
Rovira i Virgili University of Tarragona\\ Av. Pa\"{\i}sos Catalans
26, 43007 Tarragona, Spain\\ J. M.
Sigarreta\footnote{e-mail:\mbox{\tt
    josemaria.sigarreta\@@uc3m.es}}\\
{\em Departamento de Matem\'{a}ticas}\\
Universidad Carlos III de Madrid\\ Avda. de la Universidad 30, 28911
Leganés (Madrid),  Spain }

\date{}

\begin{document}

\maketitle

\begin{abstract}
In this paper we obtain several tight bounds on different types of
allian\-ce numbers of a graph, namely (global) defensive alliance
number, global offensive alliance number and  global dual alliance
number. In particular, we investigate the  relationship between
the alliance numbers of a graph and its algebraic connectivity,
its spectral radius, and its Laplacian spectral radius.
\end{abstract}



\section{Introduction}

The study of defensive alliances in graphs, together with a
variety of other kinds of allian\-ces,  was introduced by
Hedetniemi, et. al. \cite{alliancesOne}. In the referred paper was
initiated the study of the mathe\-matical properties of alliances.
In particular, several bounds on the defensive alliance number
were given. The particular case of global (strong) defensive
alliance was investigated in\cite{GlobalalliancesOne} where
several bounds on the global (strong) defensive alliance number
were obtained.

In this paper we obtain several tight bounds on different types of
allian\-ce numbers of a graph, namely (global) defensive alliance
number, global offensive alliance number and  global dual alliance
number. In particular, we investigate the  relationship between
the alliance numbers of a graph and its algebraic connectivity,
its spectral radius, and its Laplacian spectral radius.

We begin by stating some notation and terminology. In this paper
$\Gamma=(V,E)$ denotes a simple graph of  order $n$ and size $m$.
For a non-empty subset $S\subseteq V$, and any vertex $v\in V$, we
denote by $N_S(v)$ the set of neighbors $v$ has in $S$:
$$N_S(v):=\{u\in S: u\sim v\},$$
 Similarly, we denote by
$N_{V\setminus S}(v)$ the set of neighbors $v$ has in $V\setminus
S$:
$$N_{V\setminus S}(v):=\{u\in V\setminus S: u\sim v\}.$$
In this paper we will use the following obvious but useful claims:
\begin{claim} \label{claim1}
Let $\Gamma=(V,E)$ be a simple graph of size $m$. If $S\subset V,$
then \item{$$2m=\sum_{v\in S}|N_{ S}(v)|+2\sum_{v\in
S}|N_{V\setminus S}(v)|+\sum_{v\in V\setminus S}|N_{V\setminus
S}(v)|.$$}
\end{claim}

\begin{claim} \label{claim2}
Let $\Gamma=(V,E)$ be a simple graph. If $S\subset V,$ then
\item{$$\sum_{v\in S}|N_{V\setminus S}(v)|=\sum_{v\in V\setminus
S}|N_{ S}(v)|.$$}
\end{claim}

\begin{claim} \label{claim3}
Let $\Gamma=(V,E)$ be a simple graph. If $S\subset V,$ then
\item{$$\sum_{v\in S}|N_{ S}(v)|\le |S|(|S|-1).$$}
\end{claim}

\section{Defensive alliances}
 A nonempty set of vertices $S\subseteq V$ is
called a {\em defensive allian\-ce} if for every $v\in S$, $$|
N_S(v) | +1\ge  | N_{V\setminus S}(v)|.$$ In this case, by
strength of numbers,  every vertex in $S$ is {\em defended} from
possible attack by vertices in $V\setminus S$. A defensive
alliance $S$ is called {\em strong} if for every $v\in S$,
$$| N_S(v) | \ge | N_{V\setminus S}(v)|.$$ In this case
  every vertex in $S$ is {\em strongly defended}.

The {\em defensive alliance number} $a(\Gamma)$ (respectively,
{\em strong defensive allian\-ce number} $\hat{a}(\Gamma)$) is the
minimum cardinality of any defensive alliance (respectively,
strong defensive alliance) in $\Gamma$.

A particular case of alliance, called global defensive allian\-ce,
was studied in \cite{GlobalalliancesOne}. A defensive alliance $S$
is called {\em global} if it affects every vertex in $V\setminus
S$, that is, every vertex in $V\setminus S$ is adjacent to at
least one member of the alliance $S$. Note that, in this case, $S$
is a dominating set.  The {\em global defensive allian\-ce number}
$\gamma_a(\Gamma)$ (respectively, {\em global strong defensive
alliance number} $\gamma_{\hat{a}}(\Gamma)$) is the minimum
cardinality of any global defensive alliance (respectively, global
strong defensive alliance) in $\Gamma$.

\subsection{Algebraic connectivity and defensive alliances}
\label{spectralbounds}
It is well-known
that the second smallest Laplacian eigenvalue of a graph is probably
the most important information contained in the Laplacian spectrum.
This eigenvalue, frequently called {\it algebraic connectivity}, is
related to several important graph invariants and imposes reasonably
good bounds on the values of several parameters of graphs which are
very hard to compute.

 The algebraic connectivity of $\Gamma$, $\mu$, satisfies the following equality
showed by Fiedler \cite{fiedler} on weighted graphs
\begin{equation}
  \mu=2n \min \left\{ \frac{\sum_{v_i\sim v_j}(w_i-w_j)^2  }
  {\sum_{v_i\in V}\sum_{v_j\in V}(w_i-w_j)^2}: \mbox{\rm $w\neq \alpha{\bf j}$
   for   $\alpha\in \mathbb{R}$ } \right\},     \label{rfiedler}
\end{equation}
where $V=\{v_1, v_2, ..., v_n\}$, ${\bf j}=(1,1,...,1)$ and $w\in
\mathbb{R}^n$.

The following theorem  shows the relationship between the algebraic
connectivity of a graph and its (strong) defensive alliance number.

\begin{theorem} \label{cotaconnectivity}
Let $\Gamma$ be a simple graph of order $n$. Let $\mu$ be the
algebraic connectivity of $\Gamma$. The defensive alliance number
of $\Gamma$ is bounded by $$ a(\Gamma)\ge
\left\lceil\frac{n\mu}{n+\mu}\right\rceil$$ and the strong
defensive alliance number of $\Gamma$ is bounded by
$$\hat{a}(\Gamma)\ge \left\lceil\frac{n(\mu+1)}{n+\mu}\right\rceil.$$
\end{theorem}

\begin{proof}
If $S$ denotes a defensive alliance in $\Gamma$, then
\begin{equation}\label{cotavecinos}
| N_{V\setminus S}(v) | \le |S|, \quad \forall v\in S.
\end{equation}
From (\ref{rfiedler}), taking $w\in \mathbb{R}^n$ defined as
$$
w_i= \left\lbrace \begin{array}{ll} 1  & {\rm if }\quad  v_i\in S;
                            \\ 0 &  {\rm otherwise,} \end{array}
                                                \right  .$$
                                                we obtain
\begin{equation}\label{fiedlerAlliance}
\mu\le \frac{n\displaystyle\sum_{v\in S} | N_{V\setminus S}(v) |
}{|S|(n-|S|)}.
\end{equation}
Thus,  (\ref{cotavecinos}) and (\ref{fiedlerAlliance}) lead to
\begin{equation}\label{final}
\mu\le \frac{n |S|}{n-|S|}.
\end{equation}
Therefore, solving (\ref{final}) for $|S|$, and considering that
it is an integer, we obtain the bound on $a(\Gamma)$. Moreover, if
the defensive alliance $S$ is strong, then by
(\ref{fiedlerAlliance}) and Claim \ref{claim3} we obtain
\begin{equation}\label{fiedlerAlliance1}
\mu\le \frac{n\displaystyle\sum_{v\in S} | N_{S}(v) |
}{|S|(n-|S|)}\le \frac{n(|S|-1)}{n-|S|}.
\end{equation}
Hence, the result follows.
\end{proof}

The above bounds are sharp as we can check in the following
examples. It was shown in \cite{alliancesOne} that, for the
complete graph $\Gamma=K_n$, $a(K_n)=
\left\lceil\frac{n}{2}\right\rceil$ and $\hat{a}(K_n)=
\left\lceil\frac{n+1}{2}\right\rceil$. As the algebraic
connectivity of  $K_n$ is $\mu=n$, the above theorem gives the
exact value of $a(K_n)$ and $\hat{a}(K_n)$. Moreover, if $\Gamma$
is the icosahedron, then $a(\Gamma)=3$. Since in this case $n=12$
and $\mu=5-\sqrt{5}$, the above theorem gives $a(\Gamma)\ge 3$.

\begin{theorem}\label{cotafloja}
Let $\Gamma$ be a simple and connected graph of order $n$ and
maximum degree $\Delta$. Let $\mu$ be the algebraic connectivity
of $\Gamma$. The strong defensive alliance number of $\Gamma$ is
bounded by
$$\hat{a}(\Gamma)\ge \left\lceil\frac{n(\mu-\left\lfloor\frac{\Delta}{2}\right\rfloor)}{\mu}\right\rceil.$$
\end{theorem}

\begin{proof}
If $S$ denotes a strong defensive alliance in $\Gamma$, then
\begin{equation} \label{strongGrado}
|N_{V\setminus S}(v)|\le \left\lfloor\frac{deg(v)}{2}\right\rfloor
\quad \forall v\in S.
\end{equation} Thus, by
(\ref{fiedlerAlliance}) the result follows.
\end{proof}

The bound is attained, for instance, in the the following cases:
the complete graph $\Gamma=K_n$, the Petersen graph,  and the
3-cube graph.

%Note that there are cases in which the quantity
%$\mu-\left\lfloor\frac{\Delta}{2}\right\rfloor$ becomes a negative
%number. Even so,
%As the following corollary shows, there are some interesting
%particular cases of above proposition.

\subsection{Bounds on the global defensive alliance number}

The spectral radius  of a graph  is the largest eigenvalue of its
adjacency matrix. It is well-known that the spectral radius of a
graph is directly related with several parameters of the graph. The
following theorem  shows the relationship between the spectral
radius of a graph and its global (strong) defensive  alliance
number.

\begin{theorem}
Let $\Gamma$ be a simple graph of order $n$. Let $\lambda$ be the
spectral radius of $\Gamma$. The global defensive alliance number
of $\Gamma$  is bounded by
$$\gamma_{{a}}(\Gamma)\ge \left\lceil\frac{n}{\lambda+2}\right\rceil$$
and  the  global strong defensive alliance number of $\Gamma$  is
bounded by
$$\gamma_{{\hat{a}}}(\Gamma)\ge \left\lceil\frac{n}{\lambda+1}\right\rceil .$$
\end{theorem}

\begin{proof}
If $S$ denotes a  defensive alliance in $\Gamma$, then
\begin{equation}\label{cotavecinos2}
\sum_{v\in S} | N_{V\setminus S}(v)| \le \sum_{v\in S} | {N_S}(v) |+
|S|.
\end{equation}
Moreover, if the defensive  alliance $S$  is  global, we have
\begin{equation}\label{cotavecinos3}
n - |S|\le \sum_{v\in S}| N_{V\setminus S}(v)| .
\end{equation}
Thus, by (\ref{cotavecinos2}) and  (\ref{cotavecinos3}) we obtain
\begin{equation}\label{cotasuma}
n-2|S|\le \sum_{v\in S}| {N_S}(v) |.
\end{equation}
On the other hand, if ${\bf A}$ denotes the adjacency matrix if
$\Gamma$, we have
\begin{equation}\label{specRad}
\frac{\langle{\bf A}w,w\rangle}{\langle w , w\rangle} \le \lambda,
\quad \forall w\in \mathbb{R}^n\setminus \{0\}.
\end{equation}
 Thus,  taking $w$  as in the proof of Theorem \ref{cotaconnectivity}, we obtain
\begin{equation}\label{stronglerAlliance}
\sum_{v\in S}|{N_S}(v)| \le \lambda|S|.
\end{equation}
By (\ref{cotasuma}) and (\ref{stronglerAlliance}), considering
that $|S|$ is an integer, we obtain the bound on
$\gamma_{{a}}(\Gamma)$.  Moreover, if the defensive alliance $S$
is strong , then
\begin{equation}\label{cotavecinos4}
\sum_{v\in S} | N_{V\setminus S}(v)|\le \sum_{v\in S} | {N_S}(v) |.
\end{equation}
Thus, by (\ref{cotavecinos3}),  (\ref{cotavecinos4}) and
(\ref{stronglerAlliance}), we obtain $n-|S|\le \lambda |S|$. Hence,
the result follows.
\end{proof}

\begin{figure}[h]
\begin{center}
%\vspace{2cm}
%\begin{picture}(1,1)
\caption{} \label{ej1} %\vspace{-1,5cm}
\includegraphics[angle=0, width=3.5cm]{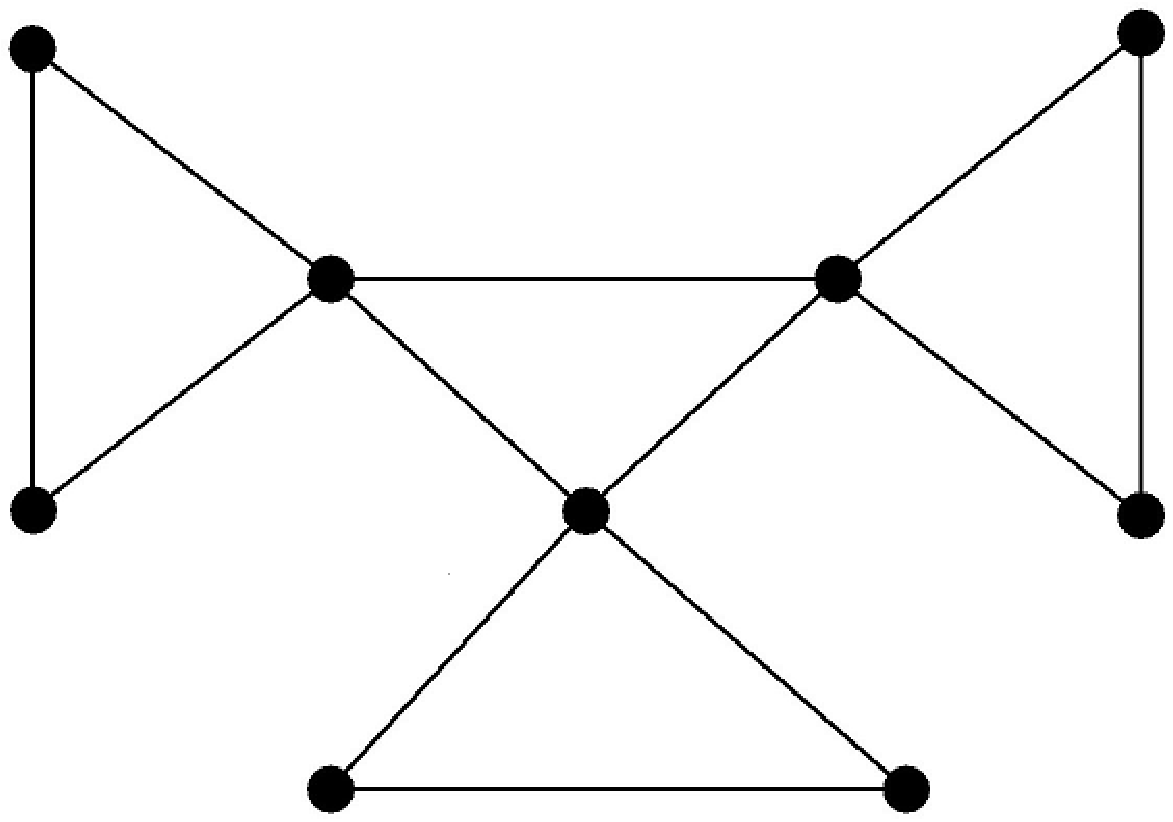}
%\epsfxsize=6cm \epsfysize=4cm \epsffile{cospectral.eps}
%\end{picture}
\end{center}
%\vspace*{-1.8cm}
\end{figure}

To show the tightness of above bounds we consider, for instance,
the graph $\Gamma=P_2\times P_3$ and the graph of Figure
\ref{ej1}. The spectral radius of $P_2\times P_3$ is
$\lambda=1+\sqrt{2}$, then we have
 $\gamma_{{a}}(\Gamma)\ge 2$.
 The spectral radius of the graph of Figure \ref{ej1} is
 $\lambda=3$, then the above theorem leads to
 $\gamma_{{\hat{a}}}(\Gamma)\ge 3$. Hence, the bounds are tight.

%XXXXXXXXXXXXXXXXX Tres triang XXXXXXXXXXXXXXXXXXXx
%\begin{figure}[h]
%\begin{center}
%\caption{ } \label{ej1}
%\includegraphics[width=0.25\textwidth]{trestriang}
%\end{center}
%\end{figure}

%\begin{figure}[h]
%\vspace*{2.cm}
%\begin{center}
%\vspace{2cm} \hspace{-6cm}
%\begin{picture}(1,1)
%\epsfxsize=6cm \epsfysize=4cm \epsffile{trestriang.eps}
%\end{picture}
%\end{center}
%\vspace*{-0.8cm}
%\end{figure}

It was shown  in \cite{GlobalalliancesOne} that if $\Gamma$ has
maximum degree $\Delta$, its global defensive alliance number  is
bounded by
\begin{equation}
\gamma_{{{a}}}(\Gamma)\ge  \frac{n}{\left\lceil
\frac{\Delta}{2}\right\rceil +1 }
\end{equation}
and  its global strong defensive alliance number is bounded by
\begin{equation}\label{cotaraiz}
\gamma_{{\hat{a}}}(\Gamma)\ge \sqrt{n}.
\end{equation}
Moreover, it was shown in  \cite{GlobalalliancesOne} that if
$\Gamma$ is bipartite, then its global  defensive alliance number is
bounded by
\begin{equation}\label{bipat}
\gamma_{{{a}}}(\Gamma)\ge \left\lceil
\frac{2n}{\Delta+3}\right\rceil.
\end{equation}
The following result  shows that the bound (\ref{bipat}) is not
restrictive to the case of bipartite graphs. Moreover, we obtain a
bound on $\gamma_{{\hat{a}}}$ that improves the bound
(\ref{cotaraiz})  in the cases of graphs of order $n$ such that
$n>\left(\left\lfloor\frac{\Delta}{2}\right\rfloor +1\right)^2$.

\begin{theorem}
Let $\Gamma$ be a simple graph of order $n$ and maximum degree
$\Delta$. The global   defensive alliance number of $\Gamma$
 is bounded by
$$\gamma_{{{a}}}(\Gamma)\ge \left\lceil \frac{2n}{\Delta+3}\right\rceil$$
and then  global strong  defensive alliance number of $\Gamma$  is
bounded by
$$\gamma_{{\hat{a}}}(\Gamma)\ge \left\lceil \frac{n}{\left\lfloor\frac{\Delta}{2}\right\rfloor+1}\right\rceil.$$
\end{theorem}

\begin{proof}
If $S$ denotes a global defensive alliance in $\Gamma$, then by
(\ref{cotavecinos3}) and (\ref{cotasuma}) we have
\begin{equation}\label{cotavecinos5}
2n-3|S|\le \sum_{v\in S} \left(| N_{V\setminus S}(v)|+ | {N_S}(v)
|\right)= \sum_{v\in S}deg(v)\le |S|\Delta.
\end{equation}
Thus, the bound on $\gamma_{{{a}}}(\Gamma)$ follows. Moreover, if
the strong defensive alliance $S$ is global, by
(\ref{cotavecinos3}) and (\ref{strongGrado}) we obtain $n\le
|S|\left(1+\left\lfloor\frac{\Delta}{2}\right\rfloor\right)$.
Hence, the bound on $\gamma_{{\hat{a}}}(\Gamma)$ follows.
\end{proof}

The tightness of the above bound of $\gamma_{a}(\Gamma)$ was showed
in \cite{GlobalalliancesOne} for the case of bipartite graphs.
Moreover, the above bound of $\gamma_{{\hat{a}}}(\Gamma)$ is
attained, for instance, in the case of the Petersen graph.

\subsection{The girth of regular graphs of small degree}

The length of a smallest cycle in a graph $\Gamma$ is called the
{\em girth} of $\Gamma$, and is denoted by $girth(\Gamma)$. It was
shown in \cite{alliancesOne} that,
\begin{enumerate}
\item[(i)]{if $\Gamma$ is regular of degree $\delta=3$ or
$\delta=4$, then $\hat{a}(\Gamma)=girth(\Gamma)$,}

\item[(ii)]{ if $\Gamma$ is 5-regular, then
${a}(\Gamma)=girth(\Gamma)$.}
\end{enumerate}

As a consequence of the previous results we obtain interesting
relations between the girth and the algebraic connectivity of
regular graphs with small degree.

\begin{theorem} \label{cubic}
Let $\Gamma$ be a simple and connected graph of order $n$. Let $\mu$
be the algebraic connectivity of $\Gamma$. Then,
\begin{itemize}
\item{
 if $\Gamma$
is 3-regular, then $girth(\Gamma)\ge
\left\lceil\frac{n(\mu-1)}{\mu}\right\rceil$;}

\item{  if $\Gamma$ is 4-regular, then $girth(\Gamma)\ge
\left\lceil\frac{n(\mu-2)}{\mu}\right\rceil;$}

\item{  if $\Gamma$ is 5-regular, then $girth(\Gamma)\ge
\left\lceil\frac{n\mu}{n+\mu}\right\rceil$.}
\end{itemize}
\end{theorem}
\begin{proof}
The results are direct consequence of (i), (ii), Theorem
\ref{cotafloja} and Theorem \ref{cotaconnectivity}.
\end{proof}

In order to show the effectiveness of above bounds we consider the
follo\-wing examples in which the bounds lead to the exact values
of the girth. If $\Gamma$ is the Petersen graph, $\delta=3$,
$n=10$ and $\mu=2$, then we have $girth(\Gamma)\ge 5$. If
$\Gamma=K_6-F$, where $F$ is a 1-factor,  $\delta=4$, $n=6$ and
$\mu=4$, then we have $girth(\Gamma)\ge 3$. If $\Gamma$ is the
 icosahedron, $\delta=5$, $n=12$ and $\mu=5-\sqrt{5}$, then we
have $girth(\Gamma)\ge 3$.
%The Petersen graph shows that the above bound is tight. In this case
%$n=10$ and $\mu=2$, Corollary \ref{cubic} leads to the exact value
%$\hat{a}(\Gamma)=5$.

\section{Offensive alliances}
The boundary of a set $S\subset V$ is defined as $$\partial
(S):=\bigcup_{v\in S}N_{V\setminus S}(v).$$ A non-empty set of
vertices $S\subseteq V$ is called {\em offensive alliance} if and
only if for every $v\in \partial (S)$,
$$| N_S(v) | \ge | N_{V\setminus S}(v)|+1.$$
An offensive alliance $S$ is called {\em strong} if for every vertex
$v\in \partial (S)$,
$$| N_S(v) | \ge | N_{V\setminus S}(v)|+2.$$
 A non-empty set of vertices $S\subseteq V$ is a {\em global offensive
alliance} if for every vertex $v\in V\setminus S$,
$$| N_S(v) | \ge | N_{V\setminus S}(v)|+1.$$
Thus, global offensive alliances are also dominating sets, and one
can define the {\em global offensive alliance number}, denoted
$\gamma_{{a}_o}(\Gamma)$, to equal the minimum cardinality of a
global offensive alliance in $\Gamma$. Analogously, $S\subseteq V$
is a {\em global strong offensive alliance} if for every vertex
$v\in V\setminus S$,
$$| N_S(v) | \ge | N_{V\setminus S}(v)|+2,$$
and the  {\em global strong offensive alliance number}, denoted
$\gamma_{\hat{a}_o}(\Gamma)$, is defined as the minimum cardinality
of a global strong offensive alliance in $\Gamma$.

\subsection{Bounds on the global offensive alliance number}
Similarly to (\ref{rfiedler}), the  Laplacian spectral radius of
$\Gamma$ (the largest Laplacian eigenvalue of $\Gamma$), $\mu_*$,
satisfies
\begin{equation}
  \mu_*=2n \max \left\{ \frac{\sum_{v_i\sim v_j}(w_i-w_j)^2  }
  {\sum_{v_i\in V}\sum_{v_j\in V}(w_i-w_j)^2}: \mbox{\rm $w\neq \alpha{\bf j}$
  for   $\alpha\in \mathbb{R}$ } \right\}.     \label{rfiedler1}
\end{equation}

The following theorem  shows the relationship between the
Laplacian spectral radius of a graph and its global (strong)
offensive alliance number.

\begin{theorem}\label{thOf}
Let $\Gamma$ be a simple graph of order $n$ and minimum degree
$\delta$. Let $\mu_*$ be the Laplacian spectral radius of
$\Gamma$. The global offensive alliance number of $\Gamma$ is
bounded by
$$\gamma_{{{a}_o}}(\Gamma)\ge \left\lceil\frac{n}{\mu_*}\left\lceil\frac{\delta +1}{2}\right\rceil\right\rceil$$
and the global strong offensive alliance number of $\Gamma$ is
bounded by
$$\gamma_{{\hat{a}_o}}(\Gamma)\ge \left\lceil\frac{n}{\mu_*}\left(\left\lceil\frac{\delta }{2}\right\rceil+1\right)\right\rceil.$$
\end{theorem}

\begin{proof}
Let $S\subseteq V$. By (\ref{rfiedler1}), taking $w\in \mathbb{R}^n$
 as in the proof of Theorem \ref{cotaconnectivity} we obtain
\begin{equation}\label{fiedlerAlliance1}
\mu_*\ge \frac{n\displaystyle\sum_{v\in V\setminus S} | N_{ S}(v) |
}{|S|(n-|S|)}.
\end{equation}
Moreover, if $S$ is a global offensive alliance  in $\Gamma$,
\begin{equation} \label{globalGrado1}
|N_{ S}(v)|\ge \left\lceil\frac{deg(v)+1}{2}\right\rceil \quad
\forall v\in V\setminus S.
\end{equation}
 Thus,
  (\ref{fiedlerAlliance1}) and (\ref{globalGrado1}) lead to
\begin{equation}\label{final1}
\mu_*\ge \frac{n}{|s|}\left\lceil\frac{ \delta  +1}{2}\right\rceil.
\end{equation}
Therefore, solving (\ref{final1}) for $|S|$, and considering that it
is an integer, we obtain the bound on $\gamma_{{{a}_o}}(\Gamma)$.
If the  global offensive alliance $S$ is strong, then
\begin{equation} \label{globalGrado2}
|N_{ S}(v)|\ge \left\lceil\frac{deg(v)}{2}\right\rceil +1\quad
\forall v\in V\setminus S.
\end{equation}
Thus,
  (\ref{fiedlerAlliance1}) and (\ref{globalGrado2}) lead to  the bound on
  $\gamma_{{\hat{a}_o}}(\Gamma)$.
\end{proof}

If $\Gamma$ is the Petersen graph, then $\mu_*=5$. Thus, Theorem
\ref{thOf} leads to $\gamma_{{{a}_o}}(\Gamma)\ge 4$ and
$\gamma_{{\hat{a}_o}}(\Gamma)\ge 6.$ Therefore, the above bounds
are tight.

\begin{theorem}
Let $\Gamma$ be a simple graph of order $n$, size $m$ and maximum
degree $\Delta$. The global offensive alliance number of $\Gamma$
is bounded by
$$\gamma_{{a_{0}}}(\Gamma)\ge \left\lceil\frac{(2n+\Delta+1)-\sqrt{(2n+\Delta+1)^{2}-8(2m+n)}}{4}\right\rceil$$
and the global strong offensive alliance number of $\Gamma$ is
bounded by
$$\gamma_{{\hat{a_{0}}}}(\Gamma)\ge \left\lceil\frac{(2n+\Delta+2)-\sqrt{(2n+\Delta+2)^{2}-16(m+n)}}{4}\right\rceil.$$
\end{theorem}

\begin{proof}
If  $S$ is a  global offensive alliance in $\Gamma=(V,E)$, then
\begin{equation}\label{ofensivaGlobal}
\sum_{v\in V\backslash S}|N_{S}(v)|\geq\sum_{v\in V\backslash
S}|N_{V\backslash S}(v)|+(n-|S|).
\end{equation}
Moreover,
 \begin{equation}\label{ofensivaglobal0}
 |S|(n-|S|)\ge \sum_{v\in V\backslash S}|N_{S}(v)|.
\end{equation}
Hence,
\begin{equation}\label{ofensivaGlobal1}
(|S|-1)(n-|S|)\ge \sum_{v\in V\backslash S}|N_{V\backslash S}(v)|.
\end{equation}
Thus,
\begin{equation}\label{d}
 (2|S|-1)(n-|S|)\ge \sum_{v\in V\backslash S}|N_{S}(v)|+
\sum_{v\in V\backslash S}|N_{V\backslash S}(v)|=\sum_{v\in
V\backslash S}deg(v).
\end{equation}
Therefore,
\begin{equation}\label{e}
(2|S|-1)(n-|S|)+\Delta|S|\ge\sum_{v\in V\backslash
S}deg(v)+\sum_{v\in S}deg(v)=2m.
\end{equation}
Thus, the bound on $\gamma_{{a_{0}}}(\Gamma)$ follows. If the
global offensive alliance $S$ is strong, then we have
\begin{equation}\label{f}
\sum_{v\in V\backslash S}|N_{S}(v)|\ge \sum_{v\in V\backslash
S}|N_{V\backslash S}(v)|+2(n-|S|).
\end{equation}
Basically the bound on $\gamma_{{\hat{a_{0}}}}(\Gamma)$ follows as
before: by replacing (\ref{ofensivaGlobal}) by (\ref{f}).
\end{proof}
The above bounds are tight as we can see, for instance, in the
case of the complete graph $\Gamma=K_{n}$ and the complete
bipartite graph $\Gamma=K_{3,6}$, for the bound on
$\gamma_{a_{0}}(\Gamma)$, and in the case of the complete
bipartite graph $\Gamma=K_{3,3}$, for the bound on
$\gamma_{\hat{a_0}}(\Gamma)$.

%As a direct consequence of Claim \ref{claim1}, Claim \ref{claim3},
%(\ref{ofensivaglobal0})  and (\ref{ofensivaGlobal1}) we deduce the
%following result.

%\begin{theorem}Let $\Gamma$ be a simple graph of order $n$ and size $m$.  The global offensive alliance
%number of $\Gamma$ is bounded by
%$$\gamma_{{a_{0}}}(\Gamma)\ge \left\lceil\frac{3n-\sqrt{9n^{2}-8n-16m}}{4}\right\rceil$$
%and the global strong offensive alliance number of $\Gamma$ is
%bounded by
%$$\gamma_{\hat{a_{0}}}(\Gamma)\ge \left\lceil\frac{3n+1-\sqrt{9n^2-10n-16m+1}}{4}\right\rceil.$$
%\end{theorem}

%XXXXXXXXXXXXXXXXX Off Cuad XXXXXXXXXXXXXXXXXXXx
%\begin{figure}[h]
%\begin{center}
%\caption{ } \label{ej5}
%\includegraphics[width=0.25\textwidth]{k3colg}
%\includegraphics[width=0.25\textwidth]{ejglobf}
%\end{center}
%\end{figure}

%The above bounds are tight: the bound on $\gamma_{a_o}(\Gamma)$ is
%attained  for the complete graph and for the left hand side graph
%of Figure \ref{ej5} and
 %the bound on $\gamma_{\hat{a_o}}(\Gamma)$ is attained
%for the right hand side graph of Figure \ref{ej5}.

\section{Dual alliances}

An alliance is called {\em dual} if it is both defensive and
offensive. The {\em global  dual alliance number} of a graph
$\Gamma$, denoted by $\gamma_{a_d}(\Gamma)$, is defined as the
minimum cardinality of any global  dual alliance in $\Gamma$. In
the case of {\em strong} alliances we denote the global  dual
alliance number by $\gamma_{\hat{a_d}}(\Gamma)$.

\subsection{Bounds on the global dual alliance number}
\begin{theorem} \label{dualSpectralradio}
Let $\Gamma$ be a simple graph of order $n$ and size $m$. Let
$\lambda$ be the spectral radius of $\Gamma$. The global dual
alliance number is of $\Gamma$ is bounded by
$$\gamma_{a_d}(\Gamma)\ge \left\lceil\frac{2m+n}{4(\lambda+1)}\right\rceil$$
and the global strong dual alliance number is of $\Gamma$ is
bounded by
$$\gamma_{\hat{a_d}}(\Gamma)\ge \left\lceil\frac{m+n}{2\lambda+1}\right\rceil.$$
\end{theorem}

\begin{proof}
Let $S$ be a global dual alliance in $\Gamma=(V,E)$.
 Since $S$
is a global offensive alliance, $S$ satisfies
(\ref{ofensivaGlobal}). Hence, by (\ref{ofensivaGlobal}) and Claim
\ref{claim1} we obtain
$$
\sum_{v\in V\setminus S}|N_S(v)|\ge \left(2m-\sum_{v\in  S}|N_
S(v)|-2\sum_{v\in S}|N_{V\setminus S}(v)|\right)+n-|S|
$$
Moreover, since the alliance $S$ is defensive, by
(\ref{cotavecinos2}) and by Claim \ref{claim2} we have
\begin{equation}\label{DualGlobal}
 4| s |+4\sum_{v\in S}|N_{S}(v)|\ge 2m+n.
\end{equation}
 Hence, by (\ref{stronglerAlliance}), the bound on
$\gamma_{a_d}(\Gamma)$ follows. On the other hand, if the global
offensive alliance $S$ is strong, then
$$
\sum_{v\in V\setminus S}|N_S(v)|\ge \sum_{v\in V\setminus
S}|N_{V\setminus S}(v)|+2(n-|S|).
$$
Hence, by Claim \ref{claim1} we have
$$
\sum_{v\in V\setminus S}|N_S(v)|\ge \left(2m-\sum_{v\in  S}|N_
S(v)|-2\sum_{v\in S}|N_{V\setminus S}(v)|\right)+2(n-|S|).
$$
and by Claim \ref{claim2} we have
$$
\sum_{v\in  S}|N_ S(v)|+3\sum_{v\in S}|N_{V\setminus S}(v)|\ge
2m+2(n-|S|).
$$
Moreover, as the strong alliance $S$ is defensive, by
(\ref{cotavecinos4}) we have
\begin{equation}\label{agfDual+claim2+defStrong}
2\sum_{v\in  S}|N_ S(v)|\ge m+n-|S|.
\end{equation}
Hence, by (\ref{stronglerAlliance}), the bound on
$\gamma_{\hat{a_d}}(\Gamma)$ follows.
\end{proof}

\begin{figure}[h]
\begin{center}
%\vspace{2cm}
%\begin{picture}(1,1)
\caption{} \label{ej4} %\vspace{-1,5cm}
\includegraphics[angle=0, width=3.5cm]{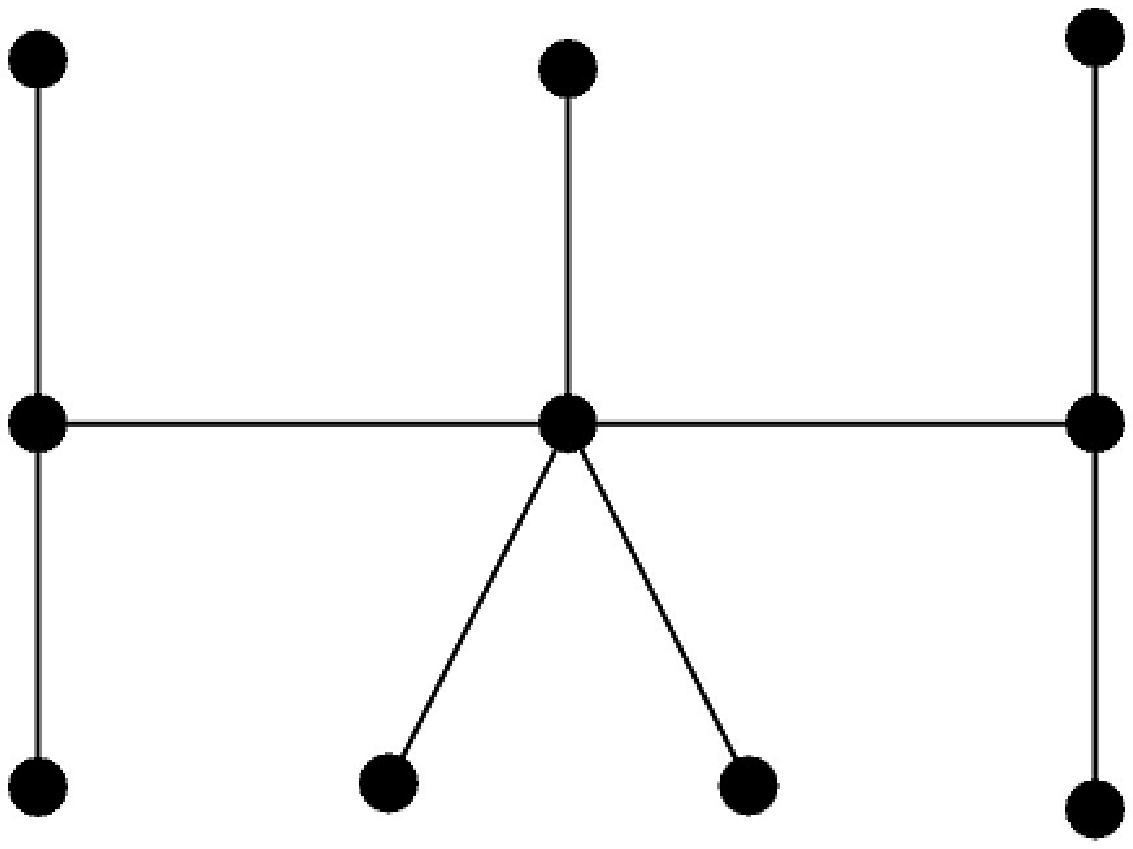}
\includegraphics[angle=0, width=3.5cm]{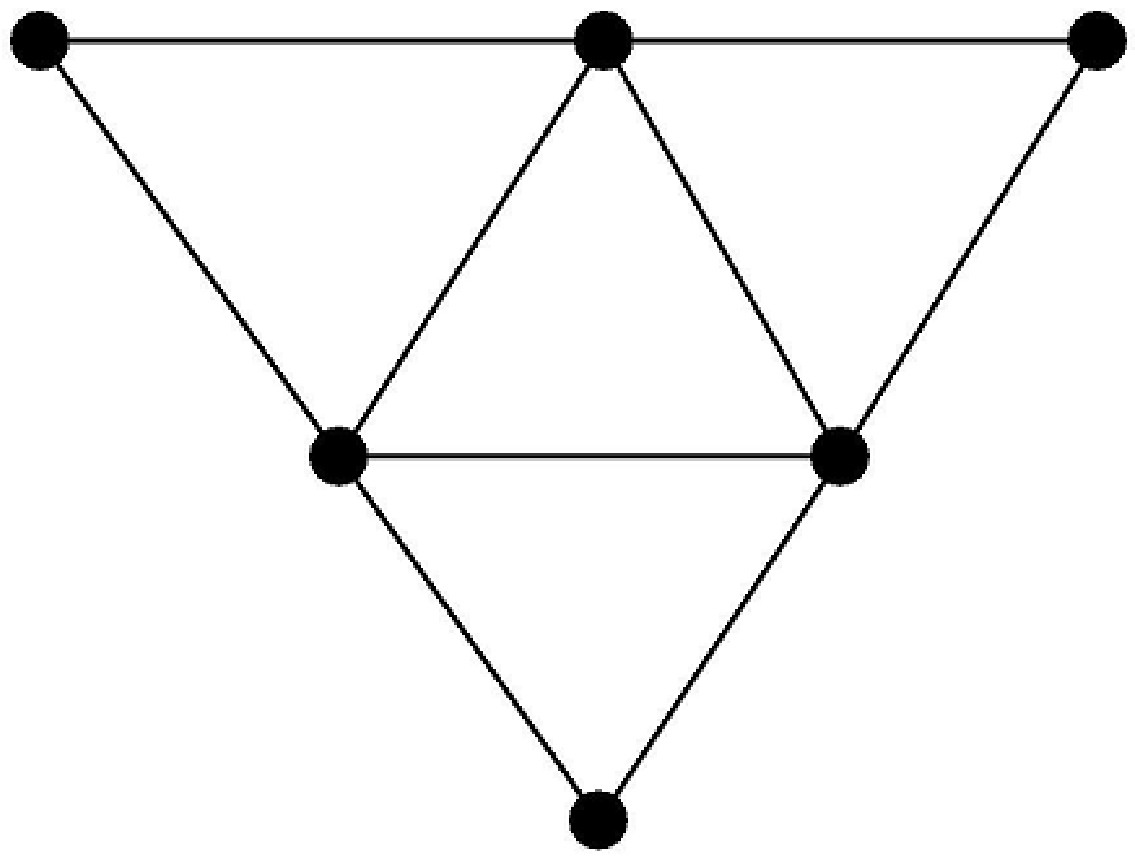}
%\epsfxsize=6cm \epsfysize=4cm \epsffile{cospectral.eps}
%\end{picture}
\end{center}
%\vspace*{-1.8cm}
\end{figure}

%XXXXXXXXXXXXXXXXX Dos Cuad XXXXXXXXXXXXXXXXXXXx
%\begin{figure}[h]
%\begin{center} \caption{ } \label{ej4}
%\includegraphics[width=0.25\textwidth]{dual}
%\includegraphics[width=0.25\textwidth]{triangulos}
%\end{center}
%\end{figure}

For the left hand side graph of Figure \ref{ej4} we have
$\lambda=\sqrt{6}$. Thus, Theorem \ref{dualSpectralradio} leads to
$\gamma_{a_d}(\Gamma)\ge 3$. Moreover, for the right hand side
graph of Figure \ref{ej4} we have $\lambda=1+\sqrt{5}$. Thus,
Theorem \ref{dualSpectralradio} leads to
$\gamma_{\hat{a_d}}(\Gamma)\ge 3$. Hence, the above bounds are
attained.

\begin{theorem}
Let $\Gamma$ be a simple graph of order $n$ and size $m$.  The
global dual alliance number is of $\Gamma$ is bounded by
$$\gamma_{a_d}(\Gamma)\ge \left\lceil\frac{\sqrt{2m+n}}{2}\right\rceil$$
and the global strong dual alliance number is of $\Gamma$ is
bounded by
$$\gamma_{\hat{a_d}}(\Gamma)\ge \left\lceil\frac{1+\sqrt{1+8(n+m)}}{4}\right\rceil.$$
\end{theorem}

\begin{proof}
Let $S$ be a global dual alliance in $\Gamma=(V,E)$.
 By  (\ref{DualGlobal}) and Claim \ref{claim3} we obtain the bound
on $\gamma_{a_d}(\Gamma)$. On the other hand, if the alliance $S$
is strong, by (\ref{agfDual+claim2+defStrong}) and Claim
\ref{claim3} we obtain the bound on $\gamma_{\hat{a_d}}(\Gamma)$.
\end{proof}

The above bounds are tight as we can see, for instance, in the
case of the complete graph $\Gamma=K_{n}$, for the bound on
$\gamma_{a_d}(\Gamma)$, and $\Gamma=K_{1}*(K_{2}\cup K_{2})$, for
the bound on $\gamma_{\hat{a_d}}(\Gamma)$, where $K_{1}*(K_{2}\cup
K_{2})$ denotes the joint of the trivial graph $K_1$ and the graph
$K_{2}\cup K_{2}$ (obtained from $K_1$ and $K_{2}\cup K_{2}$ by
joining the vertex of $K_1$ with every vertex of $K_{2}\cup
K_{2}$). Moreover, both bounds are attained in the case of the
right hand side graph of Figure \ref{ej4}.

\section{Additional observations}

By definition of global alliance, any global (defensive or
offensive) alliance is a dominating set. The {\em domination
number} of a graph $\Gamma$, denoted by $\gamma (\Gamma)$, is the
size of its smallest dominating set(s). Therefore, $\gamma_a
(\Gamma)\ge \gamma (\Gamma)$ and $\gamma_{a_o} (\Gamma)\ge \gamma
(\Gamma)$. It was shown in \cite{partition} (for the general case
of hypergraphs) that
$$\gamma (\Gamma)\ge \frac{n}{\mu_*},$$
where $\mu_*$ denotes the Laplacian spectral radius of $\Gamma$.

The reader interested in the particular case of global alliances
in planar graphs is referred to \cite{planar} for a detailed
study.

\end{document}